\documentclass[12pt,leqno]{article}
\usepackage{amsmath,amssymb}
\title{Geometry of arithmetically Gorenstein curves in ${\mathbb P}^4$}
\author{Robin Hartshorne}
\date{}

\begin{document}

\maketitle

\begin{abstract}
We characterize the postulation character of arithmetically Gorenstein
curves in ${\mathbb P}^4$.  We give conditions under which the curve can be
realized in the form $mH-K$ on some ACM surface.  Finally, we strengthen a
theorem of Watanabe by showing that any general arithmetically Gorenstein curve
in ${\mathbb P}^4$ can be obtained from a line by a series of ascending
complete-intersection biliaisons.
\end{abstract}

\setcounter{section}{-1}
\section{Introduction}
\label{sec0}

In this paper we illustrate some general results about arithmetically Gorenstein
(AG) schemes in codimension $3$ by a closer analysis of the geometry of AG
curves in ${\mathbb P}^4$.  We give a numerical criterion for when an AG curve
$Y$ in ${\mathbb P}^4$ can be obtained in the form $mH-K$ on an ACM surface $X$
whose postulation character is the expected ``first half'' of the postulation
character of $Y$.  We also give examples of AG curves that cannot be obtained in
this form on any ACM surface.  Then we prove a theorem showing that a general AG
curve with a given postulation character $\gamma$ can be obtained by ascending
CI-biliaison from a line, strengthening the result of Watanabe that any
codimension $3$ AG scheme is licci.

Section~\ref{sec1} contains a review of the structure of codimension $2$ ACM
schemes in
${\mathbb P}^N$.  Section~\ref{sec2} contains a review of known results about
codimension
$3$ AG schemes in ${\mathbb P}^N$.  Section~\ref{sec3} contains the study of AG
curves of the form $mH-K$ on an ACM surface, together with examples. 
Section~\ref{sec4} contains the theorem about ascending CI-biliaisons.

I would like to thank the University of Barcelona whose invitation to speak
there provided the impetus for writing this paper.

\section{ACM codimension 2 subschemes of ${\mathbb P}^N$}
\label{sec1}

There are in the literature many ways of recording the numerical information
associated with a subscheme of projective space.  Ellingsrud \cite{E} uses a {\em
numerical type} $(n_{ij})$ associated with a resolution of the ideal.  Gruson
and Peskine \cite{GP} use a {\em numerical character} associated with the
projection on a hyperplane.  Then there is the {\em Hilbert function} $H(n) =
h^0({\mathcal O}_X(n))$ of $X \subseteq {\mathbb P}^N$.  If $r = \dim X$, the
difference function $\partial^{r+1} H(n)$ is called the $h$-{\em vector} of $X$
\cite[\S 1.4]{M}.  Finally, there is the {\em postulation character}, also
called $\gamma$-{\em character} of \cite{MDP}.  The numerical information in
each of these is more or less equivalent, but unfortunately the terminology
varies from one place to another.  We find it most convenient, following
\cite{MDP} and \cite{N} to use the postulation character.

\bigskip
\noindent
{\bf Definition.} Let $X$ be a closed subscheme of ${\mathbb P}_k^N$ with ideal
sheaf ${\mathcal I}_X$.  Let $\varphi(n) = h^0({\mathcal O}_{{\mathbb P}^N}(n))
- h^0({\mathcal I}_X(n))$ be the postulation function of $X$.  We then define
the {\em postulation character} or $\gamma$-{\em character} of $X$ to be
$\gamma_X(n) = -\partial^N\varphi(n)$, the $N$-th difference function.

\bigskip
\noindent
{\bf Proposition 1.1} \cite{MDP}, \cite{N}.  {\em For any closed subscheme $X
\subseteq {\mathbb P}^N$, the $\gamma$-character has the following properties}
\begin{itemize}
\item[a)] {\em $\gamma(n) = 0$ for $n < 0$.}

\item[b)] {\em $\gamma(n) = -1$ for $0 \le n < s$, where $s = \min\{n \mid
h^0({\mathcal I}_X(n)) \ne 0\}$, namely, the least degree of a hypersurface
containing $X$.}

\item[c)] {\em $\gamma(s) \ge 0$.}

\item[d)] {\em $\sum_{n \in {\mathbb Z}} \gamma(n) = 0$.}
\end{itemize}

\bigskip
Note that the function $\varphi(n)$, and hence the Hilbert
polynomial of
$X$, can be recovered by numerical integration so that, for example, the degree
and genus of a curve can be expressed in terms of its $\gamma$-character.

For an ACM subscheme of ${\mathbb P}^N$ of dimension $\ge 1$ we have
$H^1({\mathcal I}_X(n)) = 0$ for all $n$, so to know the $\gamma$-character is
equivalent to knowing the Hilbert function $h^0({\mathcal O}_X(n))$.

For ACM subschemes in codimension $2$ one has precise information
about the possible $\gamma$-characters.  We call a numerical function
$\gamma(n)$ {\em admissible} if it satisfies the four conditions of $(1.1)$
for some positive integer $s$.

\bigskip
\noindent
{\bf Theorem 1.2.}  a) {\em Let $X$ be a codimension $2$} ACM {\em subscheme of
${\mathbb P}^N$.  Then its $\gamma$-character is} positive {\em in the sense
that $\gamma(n) \ge 0$ for all $n \ge s$.}

b) {\em Conversely, given an admissible numerical function $\gamma(n)$ that is
positive, as defined in} a), {\em there exists a codimension $2$} ACM {\em
subscheme $X$ in ${\mathbb P}^N$ with that $\gamma$-character.}

c) {\em If, furthermore, $X$ is integral, then $\gamma_X$ is} connected {\em in
the sense that $\{n \mid \gamma(n) > 0\}$ is a connected set of integers.}

d) {\em If $\gamma$ is a positive connected numerical function, then there
exists an integral} ACM {\em codimension $2$ subscheme $X \subseteq {\mathbb
P}^N$ with that $\gamma$-character for all $N \ge 3$.  If $N = 3$ or $4$, $X$ can
be taken to be a smooth curve or surface, respectively.}

\bigskip
\noindent
{\em Proof.} These results are due to Gruson and Peskine \cite{GP}, see also
\cite[pp.~34,111]{MDP}, \cite{M}, and \cite{N}. 

\bigskip
\noindent
{\bf Remark.} The condition ``$\gamma$ connected'' is equivalent to the condition
that the numerical character of \cite{GP} should have {\em no gaps}, and the
condition that the $h$-vector should be {\em of decreasing type}
\cite[pp.~32,97]{M}.  It is also equivalent to the condition $m(X) \ge 3$ of
Sauer (see \cite[p.~452]{N}).

\bigskip
We have also the results of Ellingsrud about the Hilbert scheme and the theorem
of Gaeta.

\bigskip
\noindent
{\bf Theorem 1.3} \cite{E}.  {\em For any positive admissible
$\gamma$-character, the set of all} ACM {\em codimension $2$ subschemes of
${\mathbb P}^N$ is a smooth, irreducible, open subset of the Hilbert scheme of
all closed subschemes of ${\mathbb P}^N$.  (There is also an explicit formula
for its dimension.)}

\bigskip
\noindent
{\bf Theorem 1.4} (Gaeta: see \cite{PS}).  {\em A codimension $2$ subscheme $X
\subseteq {\mathbb P}^N$ is} ACM {\em if and only if it is in the liaison
equivalence class of a complete intersection.}

\section{AG codimension 3 subschemes of ${\mathbb P}^N$}
\label{sec2}

Here we review the analogous results for arithmetically Gorenstein subschemes of
codimension $3$ of ${\mathbb P}^N$.  See also \cite[\S 4.3]{M} for a summary of
these results.

A closed subscheme $X \subseteq {\mathbb P}_k^N$ is {\em arithmetically
Gorenstein} (AG) if its homogeneous coordinate ring is a Gorenstein ring.  This
is equivalent to saying $X$ is ACM and its canonical sheaf $\omega_X$ is
isomorphic to ${\mathcal O}_X(m)$ for some $m \in {\mathbb Z}$.

Watanabe \cite{W} showed that the homogeneous ideal of a codimension $3$ AG
scheme is minimally generated by an odd number of elements.  His method of proof
also allows one to deduce the following result, neither stated nor proved in his
paper, but usually attributed to him.

\bigskip
\noindent
{\bf Proposition 2.1} (Watanabe).  {\em Any} AG {\em codimension $3$ subscheme
of ${\mathbb P}^N$ is in the} CI-{\em liaison class of a complete intersection.}

\bigskip
Buchsbaum and Eisenbud \cite{BE} explained the theorems of Watanabe by giving a
structure theorem for Gorenstein codimension $3$ algebras from which all further
results about AG codimension $3$ schemes are deduced.

\bigskip
\noindent
{\bf Theorem 2.2} \cite{BE}.  {\em The homogeneous ideal of any} AG {\em
codimension $3$ subscheme of ${\mathbb P}^N$ is generated by the Pfaffians of
the $(n-1) \times (n-1)$ minors of a certain skew symmetric matrix of homogeneous
polynomials, of odd rank $n$.}

\bigskip
Stanley \cite{S}, drawing on old results of Macaulay, and applying the theorem
of Buchsbaum and Eisenbud, characterized the possible $h$-vectors of AG
codimension $3$ subschemes.  Translated into the language of the
$\gamma$-character, his result is this.

\bigskip
\noindent
{\bf Proposition 2.3} \cite{S}.  {\em An admissible numerical function $\gamma$
is the $\gamma$-character of an} AG {\em codimension $3$ subscheme of ${\mathbb
P}^N$ if and only if}
\begin{itemize}
\item[a)] {\em it is} symmetric, {\em meaning there exists an integer $q$ such
that $\gamma(n) = \gamma(q-n)$ for all $n \in {\mathbb Z}$ (which implies that
$q = \max\{n \mid \gamma(n) \ne 0\}$), and}

\item[b)] {\em if we define the} $\delta$-character {\em to be the
``first half'' of $\gamma$, namely}
\[
\delta(n) = \left\{ \begin{array}{rl}
\gamma(n) &\mbox{\em for $n < \frac {q}{2}$} \\
\frac {1}{2} \gamma(n) &\mbox{\em for $n = \frac {q}{2}$} \\
0 &\mbox{\em otherwise,}
\end{array} \right.
\]
{\em then $\delta$ is a positive admissible function, as in $(1.1)$.}
\end{itemize}

\bigskip
\noindent
{\em Proof.} For a codimension $3$ subscheme of ${\mathbb P}^N$, the
$\gamma$-character is the negative of the second difference function of the
$h$-vector.  So the symmetry of the $h$-vector in Stanley's theorem
\cite[4.2]{S} is equivalent to the symmetry of $\gamma$.  For the second
condition, Stanley says the first half of the first difference function of the
$h$-vector should be an $O$-sequence.  This says it is the $h$-vector of a
codimension $2$ zero-dimensional scheme \cite[2.2]{S}.  But we know these by
$(1.2)$, namely $\delta$ should be admissible and positive.

\bigskip
\noindent
{\bf Proposition 2.4} \cite{MR}.  {\em The} AG {\em codimension $3$ subschemes of
${\mathbb P}^N$, for $N \ge 4$, are parametrized by a smooth, open subset of
the Hilbert scheme.}

\bigskip
\noindent
{\bf Proposition 2.5} \cite{D}.  {\em The} AG {\em codimension $3$ subschemes of
${\mathbb P}^N$ with a fixed Hilbert function (i.e., with a fixed
$\gamma$-character as in $(2.3)$) form an irreducible subset of the Hilbert
scheme.  (And Kleppe and Mir\'o--Roig \cite{KM} have given a formula for the
dimension of this Hilbert scheme.)}

\bigskip
For the existence of integral AG codimension $3$ schemes in ${\mathbb P}^N$,
Herzog, Trung, and Valla gave a condition in terms of the degree matrix of the
defining skew symmetric matrix of $(2.2)$.  Then De Negri and Valla translated
this condition in terms of the $h$-vector.  Combining their results and stating
them with the $\gamma$-character, we have the following.

\bigskip
\noindent
{\bf Theorem 2.6} \cite{HTV}, \cite{DV}.  a) {\em If $X$ is an integral} AG {\em
codimension $3$ subscheme of ${\mathbb P}^N$, then its $\delta$-character is
connected $(1.2c)$.}

b) {\em Conversely, if $\gamma$ is a numerical function satisfying the
conditions of $(2.3)$ with $\delta$ connected, then there exists a normal
integral} AG {\em codimension $3$ subscheme of ${\mathbb P}^N$ with that
$\gamma$-character (for $N \ge 4$).  In particular, if $N = 4$, we may take $X$
to be a smooth curve.}

\section{Arithmetically Gorenstein curves in ${\mathbb P}^4$}
\label{sec3}

Now we will look in more detail at the situation of curves in ${\mathbb P}^4$. 
If $Y$ is an AG curve in ${\mathbb P}^4$ with postulation character $\gamma_Y$,
we know from $(2.3)$ that its ``first half'' $\delta_Y$ is the postulation
character of an ACM surface $X$ in ${\mathbb P}^4$.  So our first task will be
to explore the relationship between the ACM surface $X$ and the AG curve $Y$. 
There is a well-known method of obtaining an AG curve on an ACM surface.

\bigskip
\noindent
{\bf Proposition 3.1.}  {\em Let $X$ be an} ACM {\em surface in ${\mathbb P}^4$
satisfying the additional condition $G_1$, Gorenstein in codimension $1$.  Let
$K$ be the canonical divisor, and let $Y$ be an effective divisor linearly
equivalent to $mH-K$ for some $m \in {\mathbb Z}$.  Then $Y$ is an} AG {\em curve
with
$\omega_Y \cong {\mathcal O}_Y(m)$.  If $X$ satisfies only $G_0$, then the
canonical divisor $K$ may not be defined, but there is an anticanonical divisor
$M$, and the same is true for an effective divisor $Y \sim mH+M$.}

\bigskip
\noindent
{\em Proof.}  This construction was given in \cite[5.4]{KMMNP}.  See also
\cite[4.2.8]{M}.  The extension to the case where $X$ satisfies only $G_0$ is in
\cite{H1}.

\bigskip
\noindent
{\bf Lemma 3.2.}  {\em Let $X$ be an} ACM {\em surface in ${\mathbb P}^4$ with
postulation character $\gamma_X$, and let $r = \max\{n \mid \gamma_X(n) \ne
0\}$.  Then}
\begin{itemize}
\item[1)] {\em $H^2({\mathcal O}_X(n)) = 0$ for all $n \ge r-3$.}

\item[2)] {\em ${\mathcal I}_X$ is $r$-regular.}

\item[3)] {\em ${\mathcal I}_X(r)$ is generated by global sections.}

\item[4)] {\em ${\mathcal I}_X$ has a resolution}
\[
0 \rightarrow \oplus {\mathcal O}_{{\mathbb P}^4}(-b_j) \rightarrow \oplus
{\mathcal O}_{{\mathbb P}^4}(-a_i) \rightarrow {\mathcal I}_X \rightarrow 0
\]
{\em with $\max\{b_j\} = r+1$.}
\end{itemize}

\bigskip
\noindent
{\em Proof.} Since $X$ is ACM, we have $H^1({\mathcal O}_X(n)) = 0$ for all
$n$.  Hence the Euler characteristic $\chi({\mathcal O}_X(n)) = h^0({\mathcal
O}_X(n)) + h^2({\mathcal O}_X(n))$, and this is equal to the Hilbert polynomial
of $X$.  When we take difference functions of $h^0({\mathcal O}_(n))$, the
third and fourth differences will be $0$ if and only if the corresponding shift
of $h^0({\mathcal O}_X(n))$ is equal to the polynomial $\chi({\mathcal
O}_X(n))$.  We conclude that $\gamma_X(n) = 0$ for $n \ge r+1$ is equivalent to
$h^2({\mathcal O}_X(n)) = 0$ for $n \ge r-3$.

Since ${\mathcal I}_X$ has no $H^1$ or $H^2$, because of $X$ being ACM, and
$H^3({\mathcal I}_X(n)) \cong H^2({\mathcal O}_X(n))$, we find ${\mathcal I}_X$
is $r$-regular.  This implies ${\mathcal I}(r)$ generated by global sections,
by the theorem of Castelnuovo--Mumford.

Finally, take a minimal resolution of $I_X$ over the homogeneous coordinate
ring, and sheafify.  This gives an exact sequence of cohomology
\[
0 \rightarrow H^3({\mathcal I}_X(n)) \rightarrow \oplus H^4({\mathcal
O}_{{\mathbb P}^4}(n-b_j)) \stackrel{\alpha}{\rightarrow} \oplus H^4({\mathcal
O}_{{\mathbb P}^4}(n-a_i)).
\]
Because of the minimality of the resolution, $\max\{b_j\} > \max\{a_i\}$. 
Hence the largest $n$ for which $H^3({\mathcal I}_X(n)) \ne 0$ is equal to the
largest $n$ for which some $H^4({\mathcal O}_{{\mathbb P}^4}(n-b_j)) \ne 0$.  We
conclude $n \ge r-3$ if and only if $n-b_j > 5$ for all $j$, and hence
$\max\{b_j\} = r+1$.

\bigskip
\noindent
{\bf Lemma 3.3.}  {\em Let $X$ be a locally complete intersection surface in
${\mathbb P}^4$, with ideal sheaf ${\mathcal I}$.  Then $\Lambda^2({\mathcal
I}/{\mathcal I}^2) \cong \omega_X^{\vee}(-5)$.}

\bigskip
\noindent
{\em Proof.} \cite[III.7.11]{AG}.

\bigskip
\noindent
{\bf Lemma 3.4.} {\em Let $Y$ be an} AG {\em curve in ${\mathbb P}^4$ with
$\omega_Y \cong {\mathcal O}_Y(m)$.  Then for all $n \in {\mathbb Z}$,
$\gamma_Y(n) = \gamma_Y(m+4-n)$.  In other words $\gamma_Y$ is symmetric with the
$q$ of $(2.3a)$ equal to $m+4$.}

\bigskip
\noindent
{\em Proof.} By duality on $Y$ we have $H^2({\mathcal O}_Y(n))$ dual to
$H^0({\mathcal O}_Y(m-n))$ for all $n$.  Hence by Riemann--Roch,
\[
h^0({\mathcal O}_Y(n)) = dn + 1 - g + h^0({\mathcal O}_Y(m-n)).
\]
When we take the fourth difference function, this gives
\[
\gamma_Y(n) = \gamma_Y(m+4-n).
\]
Hence $\gamma_Y$ is symmetric with $q = m+4$.

\bigskip
Now we can state our main result about the AG curves of the form $mH-K$ on an
ACM surface $X$.

\bigskip
\noindent
{\bf Theorem 3.5.}  {\em Let $X$ be an} ACM {\em surface in ${\mathbb P}^4$ that
is locally complete intersection, and let $r = \max\{n \mid \gamma_X(n) \ne
0\}$.}
\begin{itemize}
\item[a)] {\em If $m \ge 2r-5$, the linear system $|mH-K|$ is effective and
without base points.}

\item[b)] {\em If $m \ge 2r-4$, then $mH-K$ is very ample, and for any $Y \in
|mH-K|$, with postulation character $\gamma_Y$, its first half $\delta_Y$ is
equal to
$\gamma_X$.}

\item[c)] {\em If $m \ge 2r-2$, the set of curves $Y \in |mH-K|$ as $X$ also
varies in its family, forms an open subset of the Hilbert scheme if} AG {\em
curves with postulation character $\gamma_Y$.}
\end{itemize}

\bigskip
\noindent
{\em Proof.} We have assumed that $X$ is local complete intersection so that
$\omega_X$ will be an invertible sheaf and $K$ a Cartier divisor.
\begin{itemize}
\item[a)] By $(3.2)$, ${\mathcal I}_X(r)$ is generated by global sections.  Hence
the same is true of $({\mathcal I}/{\mathcal I}^2)(r)$ and also of
$\Lambda^2({\mathcal I}/{\mathcal I}^2(r))$.  Using $(3.3)$ this shows that
$\omega_X^{\vee}(2r-5)$ is generated by global sections, so the corresponding
linear system is effective and without base points.

\item[b)] It follows \cite[II. Ex. 7.5]{AG} for $m \ge 2r-4$, that $mH-K$
will be very ample.  To show that $\delta_Y = \gamma_X$, we proceed as follows. 
First suppose $m \ge 2r-3$.  I claim that $\gamma_X(n) = \gamma_Y(n)$ for $n \le
r$.  There is an exact sequence
\[
0 \rightarrow H^0({\mathcal I}_{Y,X}(n)) \rightarrow H^0({\mathcal O}_X(n))
\rightarrow H^0({\mathcal O}_Y(n)) \rightarrow H^1({\mathcal I}_{Y,X}(n))
\rightarrow \dots.
\]
We will show that the two outside terms are zero, hence the middle ones
isomorphic.  Since $Y \sim mH-K$, ${\mathcal I}_{Y,X}(n) \cong \omega_X(n-m)$. 
By duality on $X$, $h^1({\mathcal I}_{Y,X}(n)) = h^1({\mathcal O}_X(m-n)) = 0$,
since $X$ is ACM.  Also $h^0({\mathcal I}_{Y,X}(n)) = h^2({\mathcal
O}_X(m-n))$.  Now our assumptions $m \ge 2r-3$ and $n \le r$ imply $m-n \ge
r-3$, so the $h^2$ is $0$ by $(3.2)$.  Thus $h^0({\mathcal O}_X(n)) =
h^0({\mathcal O}_Y(n))$ and taking difference function $\gamma_X(n) =
\gamma_Y(n)$ for $n \le r$.

Now $q=m+4$ by $(3.4)$, hence $q \le 2r+1$, so $r < \frac {q}{2}$.  So we see
that the entire non-zero portion of $\gamma_X$ is equal to the portion of
$\delta_Y$ for $n \le r$.  Since both are positive admissible characters, they
are equal.

The same argument works for $m \ge 2r-4$ except for the case $m=2r-4$ and
$n=r$.  In this case $r = \frac {q}{2}$, and since the characters $\gamma_X$ and
$\delta_Y$ are equal for $n < r$, and $0$ for $n > r$, it follows that they are
also equal for $n=r$ by $(1.1d)$.

\item[c)] We follow the method of \cite[3.3]{H2}.  The proof given there already
shows the desired result for $m \gg 0$.  Note that in \cite{H2} the surface $X$
is supposed smooth, but the same holds for $X$ a locally complete intersection,
in which case $K$ will be a Cartier divisor.

The first step is that, since $X$ is ACM, the dimension of the linear system
$|Y|$ on $X$ is equal to $h^0(Y,{\mathcal N}_{Y/X})$.  The second step is by the
theorem of Ellingsrud $(1.3)$ to note that the dimension of the family of ACM
surfaces containing $X$ is equal to $h^0({\mathcal N}_{X/{\mathbb P}^4})$.  The
third step is to notice that for $m \ge 2r-3$, each such curve $Y$ is contained
in a unique such surface $X$.  Indeed, we have seen in the proof of b) above
that for $m \ge 2r-3$ we have $\gamma_X(n) = \gamma_Y(n)$ for $n \le r$, and
hence $H^0({\mathcal I}_X(n)) \rightarrow H^0({\mathcal I}_Y(n))$ is an
isomorphism for $n \le r$.  Since ${\mathcal I}_X$ is generated in degrees $\le
r$ $(3.2)$, we see that the homogeneous ideal of $X$ is uniquely determined by
$Y$.  It follows now that the dimension of the family of curves $Y \sim mH-K$ as
$X$ and $Y$ vary is equal to $h^0(Y,{\mathcal N}_{Y/X}) + h^0(X,{\mathcal
N}_{X/{\mathbb P}^4})$.

The next step is to show that for $m \ge 2r-2$ we have $h^0(X,{\mathcal
N}_{X/{\mathbb P}^4}) = h^0(Y,{\mathcal N}_{X/{\mathbb P}^4} \otimes {\mathcal
O}_Y)$.  To prove this, we need to show the vanishing of $h^i(X,{\mathcal
N}_{X/{\mathbb P}^4}(-Y))$ for $i = 0,1$.  Since $Y \sim mH-K$, these are equal,
using duality on $X$, to $h^i(X,{\mathcal I}_X/{\mathcal I}_X^2(m))$ for $i =
1,2$.  

Take a minimal resolution of ${\mathcal I}_X$, as in $(3.2)$,
\[
0 \rightarrow \oplus {\mathcal O}_{{\mathbb P}^4}(-b_j) \rightarrow \oplus
{\mathcal O}_{{\mathbb P}^4}(-a_i) \rightarrow {\mathcal I}_X \rightarrow 0.
\]
Tensoring with ${\mathcal O}_X(m)$ we get an exact sequence
\[
\oplus {\mathcal O}_X(m-b_j) \rightarrow \oplus {\mathcal O}_X(m-a_i)
\rightarrow  {\mathcal I}_X/{\mathcal I}_X^2(m) \rightarrow 0.
\]
From $(3.2)$ we know that $a_i,b_j \le r + 1$ for all $i,j$.  So the hypothesis
$m \ge 2r - 2$ implies that $m-a_i,m-b_j \ge r-3$ for all $i,j$.  So again by
$(3.2)$ it follows that $H^2({\mathcal O}_X(m-a_i)) = H^2({\mathcal O}_X(m-b_j))
= 0$ for all $i,j$.  Since $X$ is ACM, we have $H^1({\mathcal O}_X(n)) = 0$ for
all $n$, and now it follows easily that $H^i({\mathcal I}_X/{\mathcal I}_X^2(m))
= 0$ for $i = 1,2$.

Now as in the proof of \cite[3.3]{H2} we find that
\[
h^0({\mathcal N}_{Y/{\mathbb P}^4}) \le h^0({\mathcal N}_{Y/X}) + h^0({\mathcal
N}_{X/{\mathbb P}^4}).
\]
The other inequality comes from the fact that the dimension of the family of all
AG curves with the same $\gamma$-character is less than or equal to
$h^0({\mathcal N}_{Y/{\mathbb P}^4})$, by the differential study of the Hilbert
scheme.  We conclude equality, so the two families have the same dimension, and
the curves of the form $Y \sim mH-K$ with $Y,X$ varying form an open subset of
the Hilbert scheme of AG curves containing $Y$, as required.
\end{itemize}

\bigskip
\noindent
{\bf Corollary 3.6.} a) {\em For each numerical function $\gamma$ satisfying the
numerical conditions of $(2.3)$, there is an} AG {\em curve $Y$ with postulation
character $\gamma$, lying on an} ACM {\em surface $X$ with postulation character
$\gamma_X = \delta_Y$.}

b) {\em If, furthermore, the first half $\delta$ of $\gamma$ is connected, we may
take both $X$ and $Y$ to be smooth.}

c) {\em If the integers $m$ and $r$ associated with $\gamma$ satisfy $m \ge
2r-2$, then there is an open subset $V$ of the Hilbert scheme $H_{\gamma}$ of}
AG {\em curves with character $\gamma$, such that every $Y \in V$ is of the form
$Y \sim mH-K$ on an} ACM {\em surface $X$ with character $\delta$, the first half
of $\gamma$.}

\bigskip
\noindent
{\em Proof.} Note that this corollary gives an independent proof of the
existence results for AG curves quoted in $(2.3)$ and $(2.6)$.  To prove the
corollary, given $\gamma$, let $\delta$ be its first half.  Then there exists a
reduced, locally complete intersection ACM surface in ${\mathbb P}^4$ with
postulation character $\delta$ \cite[3.2]{N}.  Let $q  = \max\{n \mid \gamma(n)
\ne 0\}$ and take $m = q-4$.  If $r = \max\{n \mid \delta(n) \ne 0\}$, then $r
\le \frac {q}{2}$ by definition of the ``first half'' function, so $m \ge
2r-4$.  Then by the theorem, $mH-K$ is very ample, and any curve in the linear
system will be AG with postulation character having its first half equal to
$\gamma_X$.  By symmetry, $\gamma$ is uniquely determined by $\delta$ and $m$,
which shows that $\gamma_Y$ is equal to the $\gamma$ we started with.

If $\delta$ is connected, then $X$ can be taken to be smooth \cite[3.3]{N}, and
since $mH-K$ is very ample, we can take $Y$ to be smooth by the usual Bertini
theorem.

The last statement c) is just a reformulation of $(3.5c)$.

Now we will give some examples to illustrate Theorem $3.5$ and show that its
results are sharp.

\bigskip
\noindent
{\bf Example 3.7.} The linear system $mH-K$ may be effective even for $m <
2r-5$.  Let $X$ be the Castelnuovo surface.  Then $\gamma_X = -1\ -1\ 0\ 2$, $r =
3$, $2r-5=1$.  Take $m=0$.  Note that $-K = (3;1^8)$ in the usual notation for
divisor classes on $X$ (see, for example, \cite[3.3]{H3}).  This class is
effective and is represented by plane cubic curves $Y \subseteq X$.  They all
pass through a ninth point $Q \in X$.  Thus the linear system  $|Y|$ is not
without base points.

\bigskip
\noindent
{\bf Example 3.8.} If $m < 2r-4$, then $\delta_Y$ may not be equal to
$\gamma_X$.  One example is the plane cubic curve $Y$ of the previous $(3.7)$. 
In this case $\gamma_Y = -1\ 1\ 0\ 1\ -1$, and $\delta_Y = -1\ 1$ corresponding to
a plane
$H$.  This is not surprising since the minimal degree surface containing $Y$ is
a plane.

For a more interesting example let $X$ be the Del Pezzo surface.  Then $\gamma_X
= -1\ -1\ 1\ 1$, $r = 3$, $2r-4= 2$.  Take $m = 1$.  Then $Y \sim H-K = 2H$,
and
$Y$ is the complete intersection of three quadric hypersurfaces in ${\mathbb
P}^4$.  It is the canonical curve of degree $8$ and genus $5$, and $\gamma_Y =
-1\ -1\ 2\ 2\ -1\ -1$, so $\delta_Y = -1\ -1\ 2$, the $\gamma$-character of a
cubic scroll.  Our curve $Y$ is on its surface $X$ of minimal degree, yet its
$\delta_Y$ belongs to a surface of lower degree.  There are curves $Y' \sim H-K$
on the cubic scroll $X'$, also canonical curves $(8,5)$, but these form a proper
subfamily of all the AG $(8,5)$ curves:  they are the canonical embeddings of
trigonal curves of genus $5$.

In this example we see that while the family of AG curves $Y$ with $\gamma_Y =
-1\ -1\ 2\ 2\ -1\ -1$ is irreducible, equal to the family of canonical curves of
genus
$5$ in ${\mathbb P}^4$, there are two types:  the general one being a complete
intersection on the Del Pezzo surface and the special one lying on a cubic
scroll.  In both cases the associated character $\delta_Y$ is that of a cubic
scroll.

\bigskip
\noindent
{\bf Example 3.9.} We saw in the proof of $(3.5c)$ that for $m \ge 2r-3$, the
curve $Y \sim mH-K$ is contained in a unique ACM surface $X$ with $\gamma_X =
\delta_Y$.  Here we show that for $m = 2r-4$, the surface $X$ may not be unique.

Recall first that if $X,X'$ are two ACM surfaces without common components,
whose union is a complete intersection $Z$ of hypersurfaces of degrees $a,b$,
then the intersection $Y = X\cap X'$ is an AG curve \cite[4.2.1]{M}.  In fact, I
claim $Y \sim (a+b-5)H + M$ on $X$ where $M$ is the anticanonical divisor. 
Indeed, since $X \cup X' = Z$ is a complete intersection, and $X$ and $X'$ have
no common components, the surface $X$ satisfies $G_0$, so we can speak of the
anticanonical divisor $M$ \cite{H1}.  From the theory of liaison it follows that
${\mathcal I}_{X',Z} = {\mathcal H}om({\mathcal O}_X,{\mathcal O}_Z)$.  But
${\mathcal I}_{X',Z} = {\mathcal I}_{Y,X}$, and ${\mathcal H}om({\mathcal
O}_X,{\mathcal O}_Z) = \omega_X \otimes \omega_Z^{\vee}$.  Thus on $X$ we have $Y
\sim (a+b-5)H+M$, since $\omega_Z = {\mathcal O}_Z(a+b-5)$.  So it follows from
$(3.1)$ that $Y$ is AG with $m = a+b-5$.

Now for our example, let $X$ again be a Castelnuovo surface.  This surface is
contained in a unique quadric hypersurface $F_2$, a cone over the nonsingular
quadric surface in ${\mathbb P}^3$.  The divisor class group of $F_2$ is
${\mathbb Z}
\oplus {\mathbb Z}$, and $X$ is in the class of bidegree $(2,3)$.  Let $X'$ be
another Castelnuovo surface of bidegree $(3,2)$.  Then $X \cup X' = Z$ is a
complete intersection of $F_2$ with a quintic hypersurface $F_5$.  Therefore $Y
= X \cap X'$ is in the class of $2H-K$, so this $Y$ has $m=2$, and is not
contained in a unique Castelnuovo surface $X$.

We observe also that $X$ passes through the singular point of $F_2$, since by
Klein's theorem it cannot be a Cartier divisor on $F_2$, and this point is none
other than the point $Q$ mentioned above in $(3.7)$.  Indeed, the plane cubic
curve of $(3.7)$ is contained in a plane $\Pi$.  This plane intersects $F_2$ in
at least a plane cubic curve, so $\Pi \subseteq F_2$, and $\Pi$ must also
contain the singular point.  As the plane cubic curve moves in a pencil, so does
$\Pi$, and the only point in common is $Q$, which must therefore be the singular
point of $F_2$.

The same argument shows that $X'$ also contains $Q$, and so all the curves $Y
\sim 2H-K$ obtained by the construction $X \cap X'$ for various $X'$ contain
this same point $Q$.  On the other hand, the linear system $|2H-K|$ is very
ample by $(3.5b)$, so we see that the curves $Y$ obtained as $X \cap X'$ for
linked Castelnuovo surfaces $X$ and $X'$ are not general among all curves in the
linear system $|2H-K|$.

\bigskip
\noindent
{\bf Example 3.10.} We give an example to show that $(3.5c)$ is sharp, namely an
example of curves with $m = 2r-3$ on an ACM surface $X$ that are not general in
their Hilbert scheme.  Let $X$ be a Castelnuovo surface, with $\gamma_X =
-1\ -1\ 0\ 2$.  Then $r=3$.  We take $m = 2r-3 = 3$, and consider curves $Y \sim
3H-K$ on $X$.  These have character $\gamma = -1\ -1\ 0\ 2\ 2\ 0\ -1\ -1$  and
have degree $18$ and genus $28$.  Each such curve $Y$ is contained in a unique
quadric hypersurface $F_2$, which is the same one that contains the Castelnuovo
surface $X$, and therefore is singular, by Klein's theorem.

On the other hand, there are AG curves $Y$ of degree $18$ and genus $28$ of the
form $Y \sim 3H-K$ on the sextic $K3$ surface $X$, which is a complete
intersection of any quadric and cubic hypersurfaces, $X = F_2 \cap F_3$.  In
this case we can take $F_2$ to be smooth, so that the unique quadric
hypersurface $F_2$ containing $Y$ is smooth, and so $Y$ is not on a Castelnuovo
surface.  Thus the family of $Y \sim 3H-K$ on Castelnuovo surfaces is special in
$H_{\gamma}$.

\bigskip
\noindent
{\bf Example 3.11.} For our last example we show that for certain $\gamma$, the
general AG curve with postulation character $\gamma$ is not of the form $mH-K$
on any ACM surface.

Take $\gamma = -1\ -1\ -1\ 6\ -1\ -1\ -1$.  Then $\delta = -1\ -1\ -1\ 3$ is the
postulation character of a Bordiga surface.  There are curves $Y \sim 2H-K =
(11;3^{10})$ on a Bordiga surface $X$.  These curves have degree $14$ and genus
$15$.  The dimension of the family of all such $Y$ on Bordiga surfaces $X$ is
less than or equal to the dimension of the linear system $|Y|$ on $X$ plus the
dimension of the family $\{X\}$ of all Bordiga surfaces.  Now $\dim_X |Y| =
h^0({\mathcal N}_{Y/X}) = h^0({\mathcal O}_Y(Y))$.  We calculate $Y^2 = 31$ from
its divisor class representation on $X$.  Since $31 > 2g_Y - 2$, the linear
system ${\mathcal O}_Y(Y)$ is nonspecial, and $h^0({\mathcal O}_Y(Y)) = 31 + 1 -
15 = 17$.  The dimension of the family $\{X\}$ is $36$, by \cite{E}.  Thus the
family of all $Y \sim 2H-K$ on Bordiga surfaces has dimension $\le 17 + 36 = 53$.

On the other hand, from the general theory of the Hilbert scheme, we know that
the dimension of any irreducible component of the Hilbert scheme of curves of
degree $d$ and genus $g$ is $\ge 5d + 1-g$.  In our case, this gives $56$.  (In
fact, the dimension is exactly $56$ \cite{KM}.)  Thus the general AG curve $Y$
with given
$\gamma$ cannot be of the form
$2H-K$ on a Bordiga surface.

It remains to show that $Y$ cannot be of the form $mH-K$ on any other ACM
surface.  If $X$ is an ACM surface of degree $\delta$ and sectional genus $\pi$,
and if $Y \sim mH-K$, then we can easily compute the degree $d$ of $Y$ as
$(Y.H)$ on the surface, and this gives $d = (m+1)\delta - 2\pi + 2$.  In our
case, since $\gamma = -1\ -1\ -1\ 6\ -1\ -1\ -1$, we have $q = 6$ and $m = 2$. 
So we must have $d = 3\delta - 2\pi + 2$.  Now looking at the possible pairs
$(\delta,\pi)$, which will be the degree and genus of a nondegenerate ACM curve
in ${\mathbb P}^3$, we see (left to reader) that $3\delta - 2\pi + 2 \le 14$
always with equality only for $(\delta,\pi) = (6,3)$.  So the only way to obtain
the curve $Y$ as $mH-K$ on an ACM surface is as $2H-K$ on a Bordiga surface.

R.~Mir\'o--Roig \cite{MR1} has verified by a dimension count that there are
similar examples of AG curves of arbitrarily high degree that cannot be obtained
in the form $mH-K$ on any ACM surface.

\bigskip
\noindent
{\bf Problem 3.12.}  For those postulation characters $\gamma$ of AG curves for
which the associated $m$ and $r$ satisfy $m = 2r-3$ or $m = 2r-4$, give a
stratification of the Hilbert scheme of AG curves with character $\gamma$,
according to the least degree of an ACM surface containing the curve, the
gonality of the abstract curve, which ones are of the form $mH-K$ on an ACM
surface, and the dimensions of the strata, so as to generalize and complete the
information illustrated in examples $(3.8)$, $(3.10)$, and $(3.11)$ above.

\section{Complete intersection biliaison}
\label{sec4}

If $C$ is a curve in ${\mathbb P}^4$, recall that a {\em complete intersection}
(CI) {\em biliaison} of $C$ is obtained by taking a complete intersection
surface $X \subseteq {\mathbb P}^4$, and taking a curve $C' \sim C + hH$ on $X$,
where $H$ is the hyperplane class, and $h \in {\mathbb Z}$.  It is {\em
ascending} if $h \ge 0$.  The equivalence relation generated by these is called
CI-{\em biliaison} and it is equivalent to even CI-liaison \cite[4.4]{GD}.

In this section we will show that a general AG curve in ${\mathbb P}^4$ is
obtained by ascending CI-biliaisons from a line.  This provides a new proof and
strengthening of Watanabe's result $(2.1)$ for general AG curves.

\bigskip
\noindent
{\bf Lemma 4.1.} {\em Let $Y$ be an} AG {\em curve in ${\mathbb P}^4$ with
postulation character $\gamma$.  Let $s = \min\{n > 0 \mid \gamma(n) \ge 0\}$,
let $q = \max\{n \mid \gamma(n) \ne 0\}$, and let $m$ be the integer for which
$\omega_Y \cong {\mathcal O}_Y(m)$.  Then}
\begin{itemize}
\item[a)] {\em $m = q-4$.}

\item[b)] {\em ${\mathcal I}_Y(q-s)$ is generated by global sections.}

\item[c)] {\em ${\mathcal I}_Y$ is $(q-1)$-regular.}
\end{itemize}

\bigskip
\noindent
{\em Proof.} Part a) we recall for memory $(3.4)$.  For parts b), c) we use the
theorem of Buchsbaum--Eisenbud \cite{BE} in the notation of \cite[\S 5,
pp.~62-63]{HTV}.  Let $S$ be the homogeneous coordinate ring of ${\mathbb
P}^4$.  Then the homogeneous ideal $I_Y$ of $Y$ has a resolution of the form
\[
0 \rightarrow S(-c) \rightarrow \oplus S(-b_i) \rightarrow \oplus S(-a_i)
\rightarrow I_Y \rightarrow 0
\]
with $i = 1,2,\dots,2r+1$ for some positive integer $r$.  Moreover, this
resolution is {\em symmetric} in the sense that if we order $a_1 \le a_2 \le
\dots \le a_{2r+1}$ and $b_1 \ge b_2 \ge \dots \ge b_{2r+1}$, then $b_i = c -
a_i$ for each $i$.  Furthermore, if we let $u_{ij} = b_i - a_j$ be the
associated degree matrix, then $u_{ij} > 0$ for $i+j = 2r+3$.

To relate this to the invariants $s$ and $q$ of the $\gamma$-character, first
note that the $a_i$ are the degrees of a minimum set of generators of $I_Y$. 
Hence $a_1 = s$, which is the least degree of a generator.  By symmetry, $b_1 =
c-s$.  Computing $\omega_Y \cong {\mathcal E}xt^3({\mathcal
O}_Y,\omega_{{\mathbb P}^4})$ using this resolution, we find $\omega_Y \cong
{\mathcal O}_Y(c-5)$.  Hence $m = c-5$ and $q = c-1$.  From the inequality
$u_{2,2r+1} > 0$ we find $b_2 > a_{2r+1} = \max\{a_i\}$.  But $b_1 = c - s \ge
b_2$, so we find $\max\{a_i\} < c-s$.  Hence $\max\{a_i\} \le q-s$, and
${\mathcal I}_Y(q-s)$ is generated by global sections.

Finally, to show that ${\mathcal I}_Y$ is $(q-1)$-regular, we use this
resolution to show $h^2({\mathcal I}_Y(q-3)) = 0$ by climbing up the resolution
and using the fact that $h^4({\mathcal O}_{\mathbb P}(-c+q-3)) = h^4({\mathcal
O}_{\mathbb P}(-4)) = 0$.

\bigskip
\noindent
{\bf Theorem 4.2.} {\em For any postulation character $\gamma$ corresponding to
an} AG {\em curve in ${\mathbb P}^4$ (as in $(2.3)$), there is a nonempty open
subset $V_{\gamma}$ of the corresponding Hilbert scheme $H_{\gamma}$ of these
curves, such that any curve $Y \in V_{\gamma}$ can be obtained by strictly
ascending} CI-{\em biliaisons from a line in ${\mathbb P}^4$.}

\bigskip
\noindent
{\em Proof.} We will prove by induction on the degree, the following slightly
more precise statement.  For each $\gamma$, there is an open set $V_{\gamma}
\subseteq H_{\gamma}$ such that for any $Y \in V_{\gamma}$
\begin{itemize}
\item[(i)] There is a complete intersection surface $X = F_s \cap F_{q-s}$
containing $Y$ that is reduced, and such that for each irreducible component
$U_i$ of $X{\backslash}\mbox{Sing } X$, the intersection $Y \cap U_i \ne
\emptyset$, and

\item[(ii)] There is an AG curve $Y' \sim Y - H$ on $X$, with postulation
character $\gamma'$, such that $Y' \in V_{\gamma'}$.
\end{itemize}

To begin with, by definition of $s$, $Y$ is contained in a hypersurface $F_s$ of
degree $s$.  Since ${\mathcal I}_Y(q-s)$ is generated by global sections, there
is a hypersurface $F_{q-s}$ containing $Y$, whose intersection with $F_s$ is a
surface $X$.  Thus every $Y \in H_{\gamma}$ is contained in a complete
intersection surface $X_{s(q-s)}$, and so the property (i) is an open condition
on $H_{\gamma}$.

We start the induction with AG curves $Y$ having $s = 1$.  These are contained
in a ${\mathbb P}^3$, so they are complete intersection curves, and for these
the theorem is immediate.  We descend by the biliaison of (ii) unless $Y$ is a
plane curve, in which case we do biliaisons in the plane containing $Y$.

So now we assume $s \ge 2$.  Suppose for a moment that $Y \subseteq X$ satisfies
condition (i).  We will show that the linear system $|Y-H|$ is nonempty and
contains an AG curve $Y'$.  We use the exact sequence of \cite[2.10]{GD},
twisted by $-H$:
\[
0 \rightarrow {\mathcal O}_X(-H) \rightarrow {\mathcal L}(Y-H) \rightarrow
\omega_Y \otimes \omega_X^{\vee}(-H) \rightarrow 0.
\]
Now $\omega_X \cong {\mathcal O}_X(s+q-s-5) = {\mathcal O}_X(q-5)$, and
$\omega_Y = {\mathcal O}_Y(m) = {\mathcal O}_Y(q-4)$, so the sheaf on the right
is just ${\mathcal O}_Y$.  Since $h^0({\mathcal O}_X(-H)) = h^1({\mathcal
O}_X(-H)) = 0$, we find $h^0({\mathcal L}(Y-H)) = h^0({\mathcal O}_Y) = 1$, so
it has a unique section $s$ whose restriction to $Y$ is $1$.  From the condition
that $Y$ meets every irreducible component of $X{\backslash}\mbox{Sing } X$, and
$X$ being reduced, we conclude that $s$ is nondegenerate, and defines an
effective divisor $Y' \sim Y - H$ \cite[2.9]{GD}.  Furthermore, since $s$
restricted to
$Y$ is $1$, we find that $Y \cap Y' = \emptyset$ as subsets of $X$.

Now consider the sequence of \cite[2.10]{GD} for $Y'$:
\[
0 \rightarrow {\mathcal O}_X \rightarrow {\mathcal L}(Y') \rightarrow
\omega_{Y'} \otimes \omega_X^{\vee} \rightarrow 0.
\]
Since $Y' \sim Y - H$ and $Y \cap Y' = \emptyset$, the sheaf ${\mathcal L}(Y')$
is invertible isomorphic to ${\mathcal L}(-H)$ on $X{\backslash}Y$, which is a
neighborhood of $Y'$.  Thus $\omega_{Y'} \otimes \omega_X^{\vee} \cong {\mathcal
O}_{Y'}(-H)$.  From this we find $\omega_{Y'} \cong {\mathcal O}_{Y'}(q-6) =
{\mathcal O}_{Y'}(m-2)$.  Hence $Y'$, which is ACM by virtue of the biliaison
from $Y$ to $Y'$, is arithmetically Gorenstein.  We can compute its
$\gamma$-character
\[
\gamma'(n) = \left\{ \begin{array}{rl}
-1 &\mbox{for $0 \le n \le s-2$} \\
\gamma(s)-1 &\mbox{for $n = s-1,q-s-1$} \\
\gamma(n+1) &\mbox{for $s \le n \le q-s-2$} \\
-1 &\mbox{for $q-s \le n \le q-2$.}
\end{array} \right.
\]

Now we explain the induction step of the proof.  Given $\gamma$ with $s \ge 2$,
define a character $\gamma'$ by the recipe just given.  By the induction
hypothesis there exists an open set $V_{\gamma'} \subseteq H_{\gamma'}$ of
curves satisfying (i) and (ii).  Let $Y'$ be such a curve, and let $Y' \subseteq
X' = F_{s'} \cap F_{q'-s'}$ satisfy (i).  Note that $q' = q-2$ and $s'$ is
either $s$ or $s-1$.  So define a surface $X = (F_{s'} + H_1) \cap (F_{q'-s'} +
H_2)$ or $X = F_{s'} \cap (F_{q'-s'} + H_1 + H_2)$, where $H_1,H_2$ are
hyperplanes in general position.  Then $X$ is a reduced complete intersection
surface of degree $s(q-s)$.  

On this surface $X$, we will show, by an argument analogous to the one above,
that a general curve $Y$ in the linear system $Y' + H$ on $X$ is an AG curve. 
First we write the sequence of \cite[2.10]{GD} for $Y'$, twisted by $H$:
\[
0 \rightarrow {\mathcal O}_X(H) \rightarrow {\mathcal L}(Y'+H) \rightarrow
\omega_{Y'} \otimes \omega_X^{\vee}(H) \rightarrow 0.
\]
Knowing that $\omega_{Y'} \cong {\mathcal O}_{Y'}(m-2)$ we see as above that the
right-hand sheaf is ${\mathcal O}_{Y'}$.  Since $X$ is a complete intersection,
$H^1({\mathcal O}_X(H)) = 0$, so the map $H^0({\mathcal L}(Y'+H)) \rightarrow
H^0({\mathcal O}_{Y'}) \rightarrow 0$ is surjective.  Since the sheaf ${\mathcal
L}(Y'+H)$ has nondegenerate sections (for example corresponding to the trivial
biliaison $Y'+H$), there exists a nondegenerate section $s \in H^0({\mathcal
L}(Y'+H))$ whose image in $H^0({\mathcal O}_{Y'})$ is $1$.  Let $Y$ be the
associated divisor.  Then $Y \cap Y' = \emptyset$, and the sequence of
\cite[2.10]{GD} for $Y$ shows as above, that $\omega_Y \cong {\mathcal
O}_Y(m)$.  Hence $Y$ is an AG curve.  Since the trivial biliaison $Y' + H$
satisfies (i), and this is an open condition, we can choose $Y$ also so that it
satisfies (i).

Thus there exists an open subset of curves $Y \in H_{\gamma}$ satisfying (i). 
Since the procedures of constructing $Y'$ from $Y$ and $Y$ from $Y'$ are
reversible, we can find an open subset $V_{\gamma} \subseteq H_{\gamma}$ of AG
curves $Y$ satisfying (i) with the associated curve $Y'$ lying in $V_{\gamma'}$.

This completes the inductive proof of (i) and (ii).  To prove the theorem, we
take a $Y \in V_{\gamma}$, and by (ii) find a $Y' \in V_{\gamma'}$ with smaller
degree.  We continue this process until either the degree is $1$ or $s = 1$,
which we have discussed above.

\bigskip
\noindent
{\bf Remark 4.3.}  If we restrict our attention to nonsingular AG curves $Y$,
the curve will lie on a nonsingular complete intersection surface $X_{s(q-s)}$,
and if $Y$ is sufficiently general, the associated curve $Y' \sim Y-H$ will
also be nonsingular.  Note that if $Y$ has connected $\delta$-character $(2.6)$
then $Y'$ also has connected $\delta$-character.  Thus we can carry out the
CI-biliaisons using only nonsingular AG curves lying on nonsingular complete
intersection surfaces.

\bigskip
\noindent
{\bf Example 4.4.} We have seen $(3.11)$ that the general AG curve $Y$ with
$\gamma = -1\ -1\ -1\ 6\ -1\ -1\ -1$, a curve of degree $14$ and genus $15$,
cannot be obtained in the form $mH-K$ on any ACM surface.  But applying our
theorem, we see that it can be obtained in the form $Y' + H$ on a complete
intersection surface $X_{3.3}$, where $Y'$ has $\gamma' = -1\ -1\ 4\ -1\ -1$. 
This is an AG curve of degree $5$ and genus $1$, which we may take to be
nonsingular when $Y$ is sufficiently general.

The curve $Y'$ in turn lies on a complete intersection $X'_{2.2}$, a Del Pezzo
surface, and $Y'' = Y' - H'$ on $X'$ is a line.  Thus $Y$ is obtained by two
ascending CI-biliaisons from a line.

\end{document}